\documentclass[12pt]{amsart}
\usepackage[dvipdf]{graphicx}
\usepackage{setspace}
\usepackage{mathrsfs}
\usepackage{amsmath,amssymb,amsthm,amsfonts,latexsym}
\usepackage{pstricks}
\usepackage[dvips]{geometry}

\newtheorem{thm}{Theorem}[section]
\newtheorem{prop}[thm]{Proposition}

\newtheorem{lem}[thm]{Lemma}

\theoremstyle{definition}
\newtheorem{conj}[thm]{Conjecture}
\newtheorem{rem}[thm]{Remark}
\newtheorem{ex}[thm]{Example}

\numberwithin{equation}{section} 

\newcommand{\bT}{\mathbf{Tab}}
\newcommand{\bR}{\mathbf R}
\newcommand{\bS}{\mathbf{S}}
\newcommand{\pr}{^{\prime}}
\newcommand{\prpr}{^{\prime\prime}}
\newcommand{\codim}{{\rm codim\,}}
\newcommand{\Rank}{{\rm Rank\,}}
\newcommand{\sh}{{\rm sh\,}}

\newcommand{\X}{\mathcal{X}}
\newcommand{\W}{\mathcal{W}}
\newcommand{\R}{\mathcal{R}}
\newcommand{\N}{\mathcal{N}}

\newcommand{\V}{\mathcal{V}}
\newcommand{\F}{\mathcal{F}}
\newcommand{\E}{\mathcal{E}}
\newcommand{\Or}{\mathcal{O}}

\newcommand{\G}{\mathfrak{g}}
\newcommand{\Sl}{\mathfrak{s}\mathfrak{l}}
\newcommand{\gh}{\mathfrak{h}}

\newcommand{\bo}{\mathfrak{b}}

\newcommand{\nil}{\mathfrak{n}}

\newcommand{\Ze}{\mathbb Z}

\newcommand{\rar}{\rightarrow}

\makeindex

\begin{document}

\title[intersections of components]{Reducibility of the intersections of components of
a Springer Fiber}

\author{A. Melnikov}
\address{Department of Mathematics, University of Haifa, Haifa
31905, Israel.} \email{melnikov@math.haifa.ac.il}

\author{N.G.J. Pagnon}
\address{    Universit\"{a}t Duisburg-Essen
    Fachbereich Mathematik, Campus Essen
    45117 Essen, Germany .} \email{ngoc.pagnon@uni-duisburg-essen.de}

\thanks{This work has been supported by the SFB/TR 45 ``Periods, moduli spaces and
 arithmetic of algebraic varieties'' and by the Marie Curie Research Training Networks ``Liegrits".}%


\begin{abstract}
The description of the intersections of components of a Springer fiber is
a very complex problem. Up to now only two cases have been described
completely. The complete picture for the hook case has been obtained by
 N. Spaltenstein and J.A. Vargas, and for two-row case  by
 F.Y.C. Fung. They have shown in particular that the intersection
of a pair  of components of a Springer fiber is either
irreducible or empty. In both cases all the components are
non-singular and the irreducibility of the intersections is  strongly related to the
non-singularity. As it has been shown in  \cite{M-P} a bijection between orbital varieties and components of the corresponding Springer fiber in ${\rm GL}_n$ 
extends to a bijection between the irreducible components of the
intersections of orbital varieties and the irreducible components of
the intersections of components of Springer fiber  preserving their
codimensions. Here we use this bijection to compute the intersections of
the irreducible components of Springer fibers for two-column case.
In this case  the components are in general singular. As we show the
intersection of two components is non-empty. The main result of
the paper is a  necessary and sufficient
condition for the intersection of two components of the Springer
fiber  to be irreducible in two-column case. The condition is purely combinatorial. 
As an application
of this characterization, we give  first examples  of pairs of components
with a reducible intersection having components of different  dimensions.
\end{abstract}

\maketitle

\section{Introduction}
\subsection{}\label{1.1}
Let $G$ denote the complex linear algebraic group ${\rm GL}_n$ with
Lie algebra ${\mathfrak{g}}=\mathfrak{gl}_n$ on which $G$ acts by
the adjoint action. For $g\in G$ and $u\in\G$ we denote this action
by $g.u:=gug^{-1}.$

\noindent We fix the standard triangular decomposition
$\G=\nil\oplus\gh\oplus\nil^-$ where $\nil$ is the subalgebra of
strictly upper-triangular matrices, $\nil^-$ is the subalgebra of
strictly lower triangular matrices and $\gh$ is the subalgebra of
diagonal matrices of $\G.$ Let $\bo:=\gh\oplus\nil$ be the standard Borel
subalgebra  so that $\nil$ is its nilpotent radical.
Let $B$ be the (Borel) subgroup of invertible upper-triangular
matrices in $G$ so that $\bo={\rm
Lie}(B)$.
\noindent The associated Weyl group  $W=\langle s_i\rangle_{i=1}^{n-1}$ where $s_i$ is a reflection w.r.t. a simple root $\alpha_i$ is identified with the symmetric group
$\bS_n$ by taking $s_i$ to be an
elementary permutation interchanging $i$ and $i+1.$
\index{$\mathfrak{g}$}\index{$\mathfrak{b}$}\index{$\mathfrak{h}$}
\index{$\mathfrak{n}$}\index{$\mathfrak{n^-}$} \index{$G$}
\index{$B$}  \index{$\bS_n$} \index{$s_i$}

Let $\F:=G/B$ \index{$\mathcal{F}$} denote the flag
manifold. Let $G{\times}^{B}\nil$ be the space obtained as the
quotient of $G\times {\nil}$ by the right action of $B$ given by
$(g,x).b:=(gb,b^{-1}.x)$ with  $g\in G$, $x\in {\nil}$ and $b\in B$.
By the Killing form we identify the space $G{\times}^{B}{\nil}$ with
the cotangent bundle of the flag manifold $ T^{*}(G/B)$. Let $g*x$
denote the class of $(g,x)$. The map $G{\times}^{B}{\nil}\rightarrow
\F\times \G, g\ast x\mapsto(gB,g.x)$ is an embedding which identifies
$G{\times}^{B}{\nil}$ with the following closed subvariety of
$\F\times\G $ (see. \cite[p.~19]{Slo1}):
$$
\mathcal{Y}:=\{(gB,x)\ |\ x\in g.{\nil}\} $$
 The map
$f:G{\times}^{B}{\nil}\rightarrow \G, g*x\mapsto g.x$ is called the
{\bf Springer resolution}. It embeds into  the following commutative
diagram:
\begin{center}
\begin{picture}(100,70)
\put(-5,50){$G{\times}^{B}{\nil}$}\put(90,50){$\F\times \G$}
\put(45,5){$\G$}\put(30,55){\vector(1,0){55}}\put(50,58){$i$} \put(8,45){\vector(1,-1){35}}\put(15,25){$\scriptstyle f$}
\put(85,45){\vector(-1,-1){35}}\put(70,25){$\scriptstyle pr_2$}
\end{picture}
\end{center}

\noindent where ${pr}_{2}:\F\times \G\rightarrow \G, (gB,x)\mapsto x$ is the
natural projection. Since $\F$ is complete and $i$ is closed embedding $f$ is proper (because $G/B$ is
complete) and its image is exactly $G.{\nil}=\mathcal{N}$, the {\bf
nilpotent variety} of $\G$ (cf. \cite{Spr1}).

Let $x$ be a nilpotent element in  ${\nil}$. By the diagram above we
have :\label{description}
$${\F}_{x}:=f^{-1}(x)=\{gB\in \F\ |\ x\in g.{\nil}\ \}= \{gB\in \F\
|\ g^{-1}.x\in {\nil}\ \}.\eqno{(*)}$$
\index{${\mathcal{F}}_{x}$} The variety ${\F}_{x}$ is called the
{\bf Springer fiber} above $x$. It has been studied by many authors.
Springer fibers arise as fibers of Springer's resolution of singularities of
the nilpotent cone in \cite{Slo1, Spr1}. In the course of
these investigations, Springer defined $\W$-module structures on the
rational homology groups $H_{*}(\F_{x},\mathbb{Q})$ on which also
the finite group $A(x)=Z_G(x)/{Z_G^o(x)}$ (where $Z_G(x)$ is a
stabilizer of $x$ and $Z_G^o(x)$ is its identity component) acts
compatibly. Recall  that $A(x)$ is trivial for $G={\rm GL}_n.$ 
For $d={\rm dim}(\F_{x})$ the $A(x)$-fixed subspace
$H_{2d}(\F_{x},\mathbb{Q})^{A(x)}$  of the top homology is known to
be irreducible as a $W$-module\cite{Spr3}.

In \cite{Kaz}, D. Kazhdan and G. Lusztig tried to understand
Springer's work connecting nilpotent classes and representations of
Weyl groups. Among  problems  posed there,  Conjecture 6.3 in
\cite{Kaz} has stimulated the research of the relation between the
Kazhdan-Lusztig basis and  Springer fibers.

Let $x\in{\nil}$ be a nilpotent element and let $\Or_x=G.x$ be
its orbit\index{$\Or_x$}. Consider $\Or_x\cap{\nil}.$ Its
irreducible components are called {\bf orbital varieties} associated
to $\Or_x.$ By Spaltenstein's construction \cite{Spa} $\Or_x\cap{\nil}$ is a translation of $\F_x$   (see Section
\ref{2.1}).

\subsection{}\label{1.2}
For $x\in\nil$ its Jordan form is completely defined by
$\lambda=(\lambda_1,\ldots \lambda_k)$ a partition of $n$ where
$\lambda_i$ is the length of $i-$th Jordan block. Arrange the
numbers in a partition $\lambda=(\lambda_1,\ldots \lambda_k)$ in the
decreasing order (that is
$\lambda_1\geqslant\lambda_2\geqslant\cdots\geqslant
\lambda_k\geqslant 1$) and write $J(x)=\lambda.$ \index{$J(x)$}Note
that the nilpotent orbit $\Or_x$ is completely defined by $J(x).$ We
set $\Or_{J(x)}:=\Or_x$ and $\sh(\Or_x):=J(x).$ \index{$\Or_{J(x)}$}
\index{$\Or_{\lambda}$} \index{$\sh(\Or_x)$}

In turn an ordered partition can be presented as a Young diagram
$D_{\lambda}$ -- an array of $k$ rows of boxes starting on the
left with the $i$-th row containing $\lambda_i$ boxes. In such a way
there is a bijection between Springer fibers (resp. nilpotent
orbits) and Young diagrams.

Fill the boxes of Young diagram $D_\lambda$ with $n$ distinct
positive integers. \ If the entries increase in rows from left to
right and in columns from top to bottom we call such an array a
Young tableau or simply a tableau of shape $\lambda.$ Let
$\bT_{\lambda}$\index{$\bT_{\lambda}$} be the set of all Young
tableaux of shape $\lambda.$ For $T\in\bT_{\lambda}$ we put
$\sh(T):=\lambda.$ \index{$\sh(T):=\lambda$}

 By Spaltenstein
(\cite{Spa0}) and Steinberg (\cite{Ste2})
for $x\in\nil$ such that $J(x)=\lambda$ there is a bijection
between the set of irreducible components of $\F_x$ (resp. orbital
varieties associated to $\Or_\lambda$) and $\bT_\lambda$ (cf.
Section \ref{2.3}). For $T\in \bT_{\lambda}$, set
$\F_T$\index{$\F_T$} to be the corresponding component of $\F_x.$
Respectively set $\V_T$\index{$\V_T$} to be the corresponding
orbital variety associated to $\Or_\lambda.$
 Moreover, as it has been established in \cite{M-P} (cf. Section \ref{2.1}) for $T,T'\in\bT_\lambda$ the
number of irreducible components and their codimensions in
$\F_T\cap\F_{T\pr}$ is equal to the number of irreducible components
and their codimensions in $\V_T\cap \V_{T\pr}.$ Thus, the study of
intersections of  irreducible components of $\F_x$ can be reduced to the study of the intersections of orbital varieties of $\Or_x\cap\nil.$

The conjecture of Kazhdan and Lusztig mentioned above is equivalent
to the irreducibility of certain characteristic varieties
\cite[Conjecture 4]{Bor}. They have been shown to be reducible in general
by Kashiwara and Saito \cite{Kas}. Nevertheless, the description of
pairwise intersections of the irreducible components of the Springer
fibers is still open.

The complete picture of the intersections of the components
have been described by J.A. Vargas for hook case in \cite{Var}
and by F.Y.C. Fung for two-row case in \cite{Fu}. Both in  hook and two-row cases, all
 the components are non-singular, all the intersections are
irreducible or empty.

In this paper we study the components of the intersection of a pair of components
 for two-column case (that is $\lambda=(2,2,\ldots)$ ). The two-column case and the hook case are two
extreme cases in the following sense: For all  nilpotent orbits
of the given rank $k$ the orbit $\lambda=(k,1,1\ldots)$ is the most
nondegenerate and the orbit $\lambda=(2,2,\ldots)$ (with dual
partition $\lambda^*=(n-k,k)$) is the most degenerate, in the
following sense $\overline\Or_{(k,1,\ldots)}\supset
\overline\Or_{\mu}\supset\Or_{(2,\ldots,2,1,\cdots)}$ for any $\mu$ such that for $x\in\Or_\mu$ one has $\Rank x=k.$ However, it seems that
the general picture must be more close to the two-column case than to the hook case,
which is too simple and beautiful.

\subsection{}
In general we have only Steinberg's construction for orbital
varieties. Via this construction orbital varieties in $\Or_x\cap\nil$ are  as complex from  geometric point
of view as irreducible components of $\F_x$.
There is, however  a nice exception: the case of orbital varieties in
$\mathfrak{gl}_n$ associated to two-column Young diagrams. In this case
each orbital variety is a union of a finite number of $B$-orbits and we
can apply \cite{Mel1} to get the full picture of intersections of
orbital varieties. In \cite{Mel1} the special so called rank matrix
is attached to a $B$-orbit of $x\in\nil.$ In the case of $x$ of nilpotent order 2 it defines the corresponding $B$-orbit
completely. Here we use the technique of these matrices to determine
the intersection of two orbital varieties of nilpotent order two. In
particular we show that the intersection of two orbital varieties
associated to an orbit of nilpotent order 2 is not empty (see
Proposition \ref{prop2}). We give the purely combinatorial and easy
to compute necessary and sufficient condition for the irreducibility
of the intersection of two orbital varieties of nilpotent order 2
and provide some examples showing that in general such intersections
are reducible and not necessary equidimensional (see examples in
Section \ref{4.5}).

 In the subsequent paper (cf. \cite{M-P1}), we
show that the intersections of codimension 1 in two-column case are
irreducible. This together with computations in low rank cases
permits us to conjecture
\begin{conj} Given $S,T\in\bT_\lambda$. If $\codim_{\F_S}\F_T\cap\F_S=1$
then $\F_T\cap\F_S$ is irreducible.
\end{conj}

\noindent Let us now give a brief outline of the contents of the
 paper.\\

$\bullet$  To make the paper as self contained as possible  we present in
Section 2  Spaltenstein's and
Steinberg's constructions and quote the connected results essential
in further analysis.

$\bullet$ In Section 3 we provide the main result of this paper,
namely,  a purely
combinatorial necessary and sufficient condition for the
intersection of two components of the Springer fiber to be
irreducible in two-column case; as an application of this characterization, we give the
first examples for which the intersections of two components of the
Springer fiber are reducible and are not of pure dimension. This is
the most technical part of the paper.

$\bullet$ In Section 4 we give some other counter-examples concerning the possible simplification of  the construction of orbital varieties
and of their intersections in codimension one.

\section{General Construction}

\subsection{}\label{2.1}
 Given $x\in \nil$ denote $G_x=\{g\in G\ | \ g^{-1}xg\in\nil\}$. Set
$f_1:G_x\rightarrow \Or_x\cap\nil$ by $f_1(g)=g.x$ and
$f_2:G_x\rightarrow \F_x$ by $f_2(g)=gB$. Define
$\pi:\F_x\rightarrow \Or_x\cap\nil,\
gB\mapsto\pi(gB):=f_1(f_2^{-1}(gB))$. By Spaltenstein
$\pi$ induces a surjection $\hat{\pi}$ from the set of irreducible
components of $\F_x$ onto the set of irreducible components of
$\Or_x\cap\nil$, moreover the fiber of this surjective map is
exactly an orbit under the action of the component group
$A(x):=Z_G(x)/Z_G^o(x)$ (cf. \cite{Spa}). He showed also that $\F_x$
and $\Or_x\cap\nil$ are equidimensional and got the following
relations:

\begin{eqnarray}
\dim(\Or_x\cap\nil)+\dim(Z_G(x))&=&\dim(\F_x)+\dim(B)\\
\dim(\Or_x\cap\nil)+\dim(\F_x)&=&\dim(\nil)\\
 \dim(\Or_x\cap\nil)&=&\frac{1}{2}\dim(\Or_x).
\end{eqnarray}

In our setting, for the case $G={\rm GL}_n$, the component is always
trivial, so $\hat{\pi}$ is actually a bijection. As an extension of
his work, we established in \cite{M-P} the following result

\begin{prop}
Let $x\in \nil$ and let $\F_1,\F_2$ be two irreducible components of
$\F_x$ and $\V_1=\pi(\F_1),\ \V_2=\pi(\F_2)$ the corresponding orbital varieties. Let $\{\E_l\}_{l=1}^t$ be the set of irreducible components of
$\F_1\cap\F_2.$ Then $\{\pi(\E_l)\}_{l=1}^t$ is exactly the set of
irreducible components of $\V_1\cap\V_2$ and
$\codim_{\F_1}(\E_l)=\codim_{\V_1}(\pi(\E_l))$.
\end{prop}

This simple proposition shows that in the case of ${\rm GL}_n$,
orbital varieties associated to $\Or_x$ are equivalent to the
irreducible components of $\F_x.$

\subsection{}\label{2.3}

The parametrization of the irreducible components of $\F_x$ in ${\rm
GL}_n$ by standard Young tableaux is as follows.

In this case $\F$ is identified with the set of complete flags
$\xi=(V_{1}\subset\cdots V_{n}=\mathbb{C}^{n})$ and $\F_{x}\cong
\{\xi=(V_{i})\in \F\ |\ x(V_{i})\subset V_{i-1}\}$.

\noindent  Given $x\in\nil$ let $J(x)=\lambda.$ By a slight abuse of
notation we will not distinguish between the partition $\lambda$ and
its Young diagram. By R. Steinberg \cite{Ste3} and N. Spaltenstein
\cite{Spa0} we have a parametrization of the irreducible components
of ${\F}_{x}$ by the set $\bT_{\lambda}$: Let $\xi=(V_{i})\in
\F_{x}$, then we get a sutured chain in the poset of Young diagrams
$${\rm St}
(\xi):=(J(x)\supset J(x_{|V_{n-1}})\supset\cdots\supset
J(x_{|V_{2}})\supset J(x_{|V_{1}}))$$\index{${\rm St} (\ )$} where
$x_{|V_{i}}$ is the nilpotent endomorphism induced by $x$ by
restriction to the subspace $V_{i}$ and $J(x_{|V_{i+1}})$ differs
from $J(x_{|V_i})$ by one corner box. It is easy to see that the data
of such a sutured chain is equivalent to give a standard Young
tableau. So we get a map ${\rm St}:\F_{x}\rightarrow \bT_{\lambda}$.
Then the collection $ \{{\rm St}^{-1}(T)\}_{T\in \bT_{\lambda}}$ is
a partition of $\F_{x}$ into smooth irreducible subvarieties of the
same dimension and $\{\overline{{\rm St}^{-1}(T)}\}_{T\in
\bT_{\lambda}}$ are the set of the irreducible components of
${\F}_{x}$ which will be denoted by $\F_T:=\overline{{\rm
St}^{-1}(T)}$ where $T\in\bT_\lambda.$ \index{${\rm St}_{1} (\ )$}
\index{$\F_T$}

 On the level of orbital varieties the construction is as follows.
For $1\leqslant i<j\leqslant n$ consider the canonical projections
$\pi_{i,j}:\nil_n\rightarrow \nil_{j-i+1}$ acting on a matrix by
deleting the first $i-1$ columns and rows and the last $n-j$ columns
and rows. For any $u\in \Or_\lambda\cap \nil$ set
$J_n(u):=J(u)=\lambda$ and $J_{n-i}(u):=J(\pi_{1,n-i}(u))$ for any
$i\ :\ 1\leqslant i\leqslant n-1.$ Exactly as in the previous
construction we get a standard Young tableau corresponding to the
sutured chain $(J_n(u)\supset\ldots\supset J_1(u)),$ therefore we
get a map ${\rm St}_1:\Or_\lambda\cap\nil\rightarrow \bT_{\lambda}.$
Again the collection $\{{\rm St}_1^{-1}(T)\}_{T\in\bT_\lambda}$ is a
partition of $\Or_\lambda\cap\nil$ into smooth irreducible
subvarieties of the same dimensions and $\{\overline{{\rm
St}_1^{-1}(T)}\cap \Or_\lambda\}_{T\in \bT_{\lambda}}$ are orbital
varieties associated to $\Or_\lambda.$
Put ${\mathcal
V}_T:=\overline{{\rm St}_1^{-1}(T)}\cap \Or_\lambda$ where
$T\in\bT_\lambda$; \index{$\pi_{i,j}$} \index{$\V_T$} in particular, $\coprod\limits_{\lambda \vdash n}\bT_\lambda$
parameterizes the set of orbital varieties contained in $\nil$.

\subsection{}\label{1.4}

A general construction for orbital varieties by R. Steinberg (cf.
\cite{Ste2}) is as follows. For $w\in \bS_n$ consider the subspace
$$\nil\cap {^w \nil}:=\bigoplus\limits_{\alpha \in {\R}^{+}\cap
{^w{\R}^{+}}}\G_{\alpha}$$ contained in $\nil$. Then $G.(\nil\cap
{^w \nil})$ is an irreducible locally closed subvariety of the
nilpotent variety $\N$. Since $\N$ is a finite union of nilpotent
orbits, it follows that there is a unique nilpotent orbit
${\Or}$ such that $\overline{G.(\nil\cap {^w
\nil})}=\overline{{\Or}}$. Moreover
$\overline{B.(\nil\cap {^w \nil})}\cap {\Or}$
is an orbital variety associated to $\Or$ and the fundamental result in Steinberg's work is that every orbital
variety can be obtained in this way \cite{Ste2}; in particular there is a surjective
map $\varphi: \bS_n\rightarrow \coprod\limits_{\lambda \vdash n}\bT_\lambda$. The preimages of this map $\mathcal{C}_{T}:=\varphi^{-1}(T)$  are called the {\bf
geometric} (or {\bf left}) {\bf cells} of $\bS_n$. \index{$\V_w$} \index{$\Or_w$}
The geometric cells are given by Robinson-Schensted correspondence, namely for $T\in \bT_\lambda$, one has $\mathcal{C}_{T}=\{ {\rm RS}(T,S)\ :\ S\in\bT_\lambda\}$, where
RS represents the Robinson-Schensted correspondence.

\section{Two-column case}

\subsection{}\label{4.1}

In this section we use intensively the results of \cite{Mel1} and we
adopt its notation.

\noindent Set \index{$\X_2$}$\X_2:=\{x\in\nil\ | \
x^{\scriptscriptstyle 2}=0\}$ to be the variety of nilpotent
upper-triangular matrices of nilpotent order 2.   Denote
\index{$\bS_n^2$}$\bS_n^2:=\{\sigma\in\bS_n\ | \
\sigma^{\scriptscriptstyle 2}=id\}$ the set of involutions of
$\bS_n$. For every $\sigma\in\bS_n^2$, set $N_\sigma$ to be the
``strictly upper-triangular part" of its corresponding permutation
matrix, that is \index{$N_\sigma$} \index{$(N_\sigma)_{i,j}$}
\begin{eqnarray}(N_\sigma)_{i,j}:=\left\{\begin{array}{ll} 1&{\rm if}\ i<j\
{\rm and}\
\sigma(i)=j;\\
0&{\rm otherwise.}\\
\end{array}\right.
\end{eqnarray}

\noindent Let $\bT_n^2$\index{${\bf Tab}_n^2$} be the set of all
Young tableaux of size $n$ with two columns. For $T\in\bT_n^2$,
write it as $T=(T_1,T_2),$ where
$T_1=\left(\begin{array}{c}t_{1,1}\cr\vdots\cr
t_{n-k,1}\cr\end{array}\right)$ is the first column of $T$ and
$T_2=\left(\begin{array}{c}t_{1,2}\cr\vdots\cr
t_{k,2}\cr\end{array}\right)$ is the second column of $T.$ And
define the following involution

\begin{eqnarray}\label{ecriture}
\sigma_{\scriptscriptstyle T}:=(i_1,j_1)\ldots (i_k,j_k)
\end{eqnarray}\index{$\sigma_T$}

\noindent where $j_s:=t_{s,2};\ i_1:=t_{1,2}-1,$ and
$i_s:=\max\{d\in T_1-\{i_1,\ldots,i_{s-1}\}\ | \ d<j_s\}$ for any
$s>1.$  For example, take
$$T=\begin{tabular}{|c|c|}
      \hline
      1 & 4 \\
      \hline
      2 & 5 \\
      \hline
      3 & 7 \\
      \hline
      6 & 8 \\
      \hline
    \end{tabular}
$$
Then $\sigma_{\scriptscriptstyle T}=(3,4)(2,5)(6,7)(1,8).$

\begin{rem}
To define $T\in\bT_n$ it is enough to know columns $T_i$ as  sets (we denote them by $\langle T_i\rangle$), or equivalently   the different column
positions $c_T(i)$ of integers $i\ :\ 1\leqslant i\leqslant n$ since
the entries increase from up to down in the columns. Thus given
$\sigma_T$ we can reconstruct $T$. Indeed, $\langle
T_2\rangle=\{j_1,\ldots,j_k\}$ and $\langle
T_1\rangle=\{i\}_{i=1}^n\setminus\langle T_2\rangle$.
\end{rem}

 One has

\begin{thm}\label{4.1spec}(\cite[2.2]{Mx2}, \cite[4.13]{Msmith})
\begin{itemize}
\item[(i)] The variety $\X_2$ is a finite union of $B$-orbits, namely
$$\X_2=\coprod\limits_{\sigma\in\bS_n^2}B.N_{\sigma}.$$

\item[(ii)]
For any $T\in \bT_n^2$, one has
$\overline\V_T=\overline{B.N}_{\sigma_T}$.
\end{itemize}
\end{thm}

\noindent The  finiteness property is particular for $\X_2.$ The
fact that each orbital variety has a dense $B$-orbit is also particular
for very few types of nilpotent orbits including orbits of nilpotent order 2 (cf.
  \cite{Msmith}). The first property
permits us to compute the intersections of any two $B$-orbit
closures in $\X_2$. The second one permits us to apply the results
to the intersections of orbital varieties of nilpotent order 2.\\

We begin with the general theory of the intersections of
$\overline{B.N}_\sigma$ for $\sigma\in \bS_n^2.$

\subsection{}\label{4.2}

In this section we prefer to use the dual partition $\lambda^*$ instead of $\lambda$ since it will be more convenient to
write it down for nilpotent orbits of nilpotent order 2. Indeed, for
$x\in \X_2$ one has $J^*(x)=(n-k,k)$\index{$J^*(x)$} where $k$ is
number of Jordan blocks of length two in $J(x).$

\begin{rem} For every element $x\in \X_2$, the integer ${\rm rk}(x)$
is exactly the number of blocks of length 2 in $J(x)$, so it defines the ${\rm GL}_n$-orbit of $x$.
\end{rem}

Any element $\sigma\in \bS_n^2$ can be written as a product of
disjoint cycles of length 2. Order elements in increasing order
inside the cycle and order cycles in increasing order according to
the first entries.  In that way we get a unique writing of every
involution. Thus, $\sigma=(i_1,j_1)(i_2,j_2)\ldots(i_k,j_k)$ where
$i_s<j_s$ for any $1\leqslant s\leqslant k$ and $i_s<i_{s+1}$ for
any $1\leqslant s<k.$ Set $L(\sigma):=k$\index{$L(\ )$} [do not
confuse this notation with the length function], and denote by
$\Or_\sigma$ the ${\rm GL}_n$-orbit of $N_\sigma$. By definition we
have $L(\sigma)={\rm rk}(N_\sigma).$ \index{$\Or_{\sigma}$}

Let us define the following number

\begin{eqnarray}r_s(\sigma):={\rm card}\{i_p<i_s\ | \ j_p<j_s\}+ {\rm
card}\{j_p\ |\ j_p<i_s\}.
\end{eqnarray}\index{$r_s(\sigma)$}

\noindent Note that the definition of $r_s(\sigma)$ is independent
of ordering cycles in increasing order according to the first
entries. However if it is ordered then $r_1(\sigma)=0$ and to
compute $r_s(\sigma)$ it is enough to check only the pairs
$(i_p,j_p)$ where $p<s.$  For example, take
$\sigma=(1,6)(3,4)(5,7).$ Then $L(\sigma)=3$ and $r_1(\sigma)=0,\
r_2(\sigma)=0,\ r_3(\sigma)=2+1=3.$

\noindent By \cite[3.1]{Mx2} one has

\begin{thm} For $\sigma=(i_1,j_1)(i_2,j_2)\ldots(i_k,j_k)\in \bS_n^2$ one has
$$\dim(B.N_\sigma)=kn-\sum\limits_{s=1}^k(j_s-i_s)-\sum\limits_{s=2}^k r_s(\sigma).$$
\end{thm}

\begin{rem}\label{rem3}
By Theorem \ref{4.1spec} (ii),the
orbits $B.N_{\sigma_{T}}$ (where $(\sh(T))^*=(n-k,k)$) are the
only $B$-orbits of maximal dimension inside the variety $\Or_{(n-k,k)^*}\cap\nil$
 and
$\dim(B.N_{\sigma_T})=k(n-k)$: Indeed any orbit $B.N_\sigma$ is
irreducible and therefore lies inside an orbital variety $\V_T$, in
particular it lies in $\overline\V_T$, so if $\dim
B.N_\sigma=\dim\V_T$ we get that $\overline {
B.N}_\sigma=\overline\V_T$ thus by Theorem \ref{4.1spec} (ii)
$\overline{B.N}_\sigma=\overline{B.N}_{\sigma_T}$ which provides
$\sigma=\sigma_T.$

In particular if $\sigma=(i_1,j_1)\ldots(i_k,j_k)$ is such that
$\dim(B.N_\sigma)=k(n-k)$, then $\sigma=\sigma_T$ where $T$ is the
tableau obtained by
$$c_T(s)=\left\{\begin{array}{ll} 2&{\rm if\ } s=j_p\ {\rm for\ some\ }p: \ 1\leqslant p\leqslant k\\
                                  1&{\rm otherwise}\\
                                  \end{array}\right.$$

\end{rem}

\subsection{}\label{4.3}

In \cite{Mel1} the combinatorial description of
$\overline{B.N}_\sigma$ (with respect to Zariski topology) for
$\sigma\in\bS_n^2$ is provided. Let us formulate this result.

Recall from Section \ref{2.3} the notion
$\pi_{i,j}:\nil_n\rightarrow \nil_{j-i+1}$ and define the rank
matrix $R_x$ of $x\in\nil$ to be

\begin{eqnarray}(R_x)_{i,j}:=\left\{ \begin{array}{ll}
0&{\rm if}\ i\geqslant j;\\
                       {\rm rk} \,(\pi_{i,j}(x))&{\rm otherwise}.\\
                       \end{array}\right.
\end{eqnarray}\index{$R_x$} \index{$(R_x)_{i,j}$}

\noindent Note that for any element $b\in B$, $\pi_{i,j}(b)$ is an
invertible upper-triangular matrix in ${\rm GL}_{j-i+1}$. Therefore
we can define an action of $B$ on $\nil_{j-i+1}$ by:
$b.y:=\pi_{i,j}(b).y$ for $y\in\nil_{j-i+1}$ and $b\in B$.

Let us first establish a result

\begin{lem}
\begin{enumerate}
\item[(i)]If $x,y\in \nil$ are in the same $B$-orbit, then they have the
same rank matrix.
\item[(ii)] The morphism $\pi_{i,j}$ is $B$-invariant.
\end{enumerate}
\end{lem}

\begin{proof}
Note that for any two upper-triangular matrices $a,b$ and for any
$i,j\ :\ 1\leqslant i<j\leqslant n$ one has
$\pi_{i,j}(ab)=\pi_{i,j}(a)\pi_{i,j}(b).$ In particular, if $a\in B$
then $\pi_{i,j}(a^{-1})=(\pi_{i,j}(a))^{-1}.$ Applying this to $x\in
\nil$ and $y$ in its $B$ orbit (that is $y=b.x$ for some $b\in B$) we
get $\pi_{i,j}(y)=\pi_{i,j}(b).\pi_{i,j}(x)$ so that the morphism
$\pi_{i,j}$ is $B$-invariant and in particular ${\rm
rk}(\pi_{i,j}(y))={\rm rk}(\pi_{i,j}(x)).$ Hence $R_x=R_y.$
\end{proof}

By this lemma we can define $R_\sigma:=R_{N_\sigma}$ as the rank
matrix associated to orbit $B.N_\sigma$. \index{$R_\sigma$}
\index{$R_{N_{\sigma}}$}
\begin{rem}\label{ones} Note that computation of $(R_{N_\sigma})_{i,j}$
is trivial -- this is exactly the number of non-zero entries in submatrix
of $1,\dots,j$ columns and $i,\ldots,n$ rows of $N_\sigma$ or in other words
the number of ones in $N_\sigma$ to the left-below of position $(i,j)$ (including
position $(i,j)$).
\end{rem}

Let $\Ze^+$ be the set of non-negative integers. Put
$\bR^2_n:=\{R_\sigma\ |\ \sigma\in \bS_n^2\}.$ By \cite[3.1,
3.3]{Mel1} one has\index{$\bR^2_n$}

\begin{prop} \label{ens}
$R=(R_{i,j})\in M_{n\times n}(\Ze^+)$  belongs to
$\bR^2_n$ if and only if it satisfies
\begin{itemize}
\item[(i)] $R_{i,j}=0$ if $i\geqslant j;$
\item[(ii)] For $i<j$ one has $R_{i+1,j}\leqslant R_{i,j}\leqslant R_{i+1,j}+1$
            and $R_{i,j-1}\leqslant R_{i,j}\leqslant R_{i,j-1}+1;$
\item[(iii)] If $R_{i,j}=R_{i+1,j}+1=R_{i,j-1}+1=R_{i+1,j-1}+1$ then
\begin{itemize}

\item[(a)] $R_{i,k}=R_{i+1,k}$ for any $k<j$ and
$R_{i,k}=R_{i+1,k}+1$ for any $k\geqslant j;$
\item[(b)] $R_{k,j}=R_{k,j-1}$ for any $k>i$ and
$R_{k,j}=R_{k,j-1}+1$ for any $k\leqslant i;$
\item[(c)] $R_{j,k}=R_{j+1,k}$ and $R_{k,i}=R_{k,i-1}$
for any $k\ :\ 1\leqslant k\leqslant n.$
\end{itemize}
\end{itemize}
\end{prop}

\noindent Fix $\sigma \in\bR^2_n$, then the conditions (i) and (ii)
 are obvious from Remark \ref{ones}, and the conditions
(iii) appears exactly for the coordinates $(i,j)$ in the matrix when
$j=\sigma(i)$, with $i<j$; we draw the following picture (see Figure
\ref{Fig4} below) to help the reader to visualize the constraints
(a), (b), (c) of (iii), with the following rule: the integers which
are inside a same white polygon, are equal, and the integers in a
same gray rectangle differ by one.

\begin{figure}[htbp]
  \centering
\psset{unit=0.5cm}
\begin{tabular}{cc}
  \begin{pspicture}(-6,-9)(7,6)

\pspolygon(0.05,-0.05)(1,-0.05)(1,-1)(1.95,-1)(1.95,-1.95)(0.05,-1.95)
\uput{0.45}[ur](0.7,-1){$+1$}

\pspolygon(-1.95,-0.05)(-1.05,-0.05)(-1.05,-1.95)(-1.95,-1.95)
\pspolygon(-0.95,-0.05)(-0.05,-0.05)(-0.05,-1.95)(-0.95,-1.95)
\psline[linestyle=dashed,linearc=1pt]{<->}(-2.5,-2)(-2.75,-1.5)(-2.75,-0.5)(-2.5,0)
\uput{0.45}[ur](-4,-1.5){$=$}
\psline[linestyle=dotted,dotsep=1pt](-3.90,-0.5)(-2.15,-0.5)
\psline[linestyle=dotted,dotsep=1pt](-3.90,-1.5)(-2.15,-1.5)

\pspolygon(-5.95,-0.05)(-4.05,-0.05)(-4.05,-1.95)(-5.95,-1.95)

\pspolygon(-5.95,-2.05)(-4.05,-2.05)(-4.05,-2.95)(-5.95,-2.95)
\pspolygon(-5.95,0.05)(-4.05,0.05)(-4.05,0.95)(-5.95,0.95)

\psline[linestyle=dotted,dotsep=1pt](-6,-0.5)(-7,-0.5)
\psline[linestyle=dotted,dotsep=1pt](-6,-1.5)(-7,-1.5)

\psline[linestyle=dotted,dotsep=1pt](-5.5,1.1)(-5.5,4.95)
\psline[linestyle=dotted,dotsep=1pt](-4.5,1.1)(-4.5,4.95)

\psline[linestyle=dotted,dotsep=1pt](-5.5,-3.1)(-5.5,-5)
\psline[linestyle=dotted,dotsep=1pt](-4.5,-3.1)(-4.5,-5)

\pspolygon(0.05,-2.05)(1.95,-2.05)(1.95,-2.95)(0.05,-2.95)
\pspolygon(0.05,-3.05)(1.95,-3.05)(1.95,-3.95)(0.05,-3.95)
\psline[linestyle=dashed,linearc=1pt]{<->}(0,-4.5)(0.5,-4.75)(1.5,-4.75)(2,-4.5)
\uput{0.45}[ur](0.4,-5.5){$=$}
\psline[linestyle=dotted,dotsep=1pt](0.5,-4.1)(0.5,-6)
\psline[linestyle=dotted,dotsep=1pt](1.5,-4.1)(1.5,-6)

\pspolygon(0.05,-6.05)(1.95,-6.05)(1.95,-7.95)(0.05,-7.95)

\psline[linestyle=dotted,dotsep=1pt](0.5,-8.1)(0.5,-9)
\psline[linestyle=dotted,dotsep=1pt](1.5,-8.1)(1.5,-9)

\pspolygon(2.05,-6.05)(2.95,-6.05)(2.95,-7.95)(2.05,-7.95)
\pspolygon(-0.05,-6.05)(-0.95,-6.05)(-0.95,-7.95)(-0.05,-7.95)

\psline[linestyle=dotted,dotsep=1pt](3.1,-6.5)(6.95,-6.5)
\psline[linestyle=dotted,dotsep=1pt](3.1,-7.5)(6.95,-7.5)

\psline[linestyle=dotted,dotsep=1pt](-1.1,-6.5)(-3,-6.5)
\psline[linestyle=dotted,dotsep=1pt](-1.1,-7.5)(-3,-7.5)

\pspolygon[fillstyle=solid,fillcolor=lightgray](2.05,0)(2.95,0)(2.95,-1.95)(2.05,-1.95)
\pspolygon[fillstyle=solid,fillcolor=lightgray](3.05,0)(3.95,0)(3.95,-1.95)(3.05,-1.95)
\psline{-}(2.05,-1)(2.2,-1)\psline{-}(2.45,-1)(2.60,-1)
\psline{-}(2.8,-1)(2.95,-1)
\psline{-}(3.05,-1)(3.2,-1)\psline{-}(3.45,-1)(3.60,-1)
\psline{-}(3.8,-1)(3.95,-1)
\psline[linestyle=dotted,dotsep=1pt](4.1,-0.5)(6,-0.5)
\psline[linestyle=dotted,dotsep=1pt](4.1,-1.5)(6,-1.5)
\psline[linestyle=dashed,linearc=1pt]{->}(4.5,-2)(4.75,-1.5)(4.75,-0.5)(4.5,0)
\uput{0.45}[ur](4.5,-1.5){\tiny{$+1$}}
\pspolygon[fillstyle=solid,fillcolor=lightgray](6.05,0)(6.95,0)(6.95,-1.95)(6.05,-1.95)
\psline{-}(6.05,-1)(6.2,-1)\psline{-}(6.45,-1)(6.60,-1)
\psline{-}(6.8,-1)(6.95,-1)

\uput{0.45}[ur](7,-1){\tiny{$i$}} \uput{0.45}[ur](7,-7){\tiny{$j$}}
\uput{0.45}[ur](7,-8){\tiny{$j+1$}}

\pspolygon[fillstyle=solid,fillcolor=lightgray](0.05,0.05)(1.95,0.05)(1.95,0.95)(0.05,0.95)
\pspolygon[fillstyle=solid,fillcolor=lightgray](0.05,1.05)(1.95,1.05)(1.95,1.95)(0.05,1.95)
\psline{-}(1,0.05)(1,0.30)
\psline{-}(1,0.45)(1,0.60)
\psline{-}(1,0.8)(1,0.95)
\psline{-}(1,1.05)(1,1.30)
\psline{-}(1,1.45)(1,1.60)
\psline{-}(1,1.8)(1,1.95)

\psline[linestyle=dotted,dotsep=1pt](0.5,2.25)(0.5,4.95)
\psline[linestyle=dotted,dotsep=1pt](1.5,2.25)(1.5,4.95)
\psline[linestyle=dashed,linearc=1pt]{->}(0,2.5)(0.5,2.75)(1.5,2.75)(2,2.5)
\uput{0.45}[ur](0.25,2.75){\tiny{$+1$}}
\pspolygon[fillstyle=solid,fillcolor=lightgray](0.05,4.05)(1.95,4.05)(1.95,4.95)(0.05,4.95)
\psline{-}(1,4.05)(1,4.30)\psline{-}(1,4.45)(1,4.60)
\psline{-}(1,4.8)(1,4.95)

\uput{0.45}[ur](1,5){\tiny{$j=\sigma(i)$}}
\uput{0.45}[ur](-4.9,5){\tiny{$i$}}
\uput{0.45}[ur](-6.5,5){\tiny{$i-1$}}

\psline{-}(-7,4.95)(6.95,4.95)(6.95,-8.95)

\end{pspicture}
\end{tabular}\caption{} \label{Fig4}
\end{figure}

The first part of (c) can be explained in the following: since the
integer $j$ appears already in the second entry of the cycle
$(i,j)$, so it can not appear again in any other cycle; therefore in
the matrix $N_{\sigma}$, the integers of the $j^{th}$ row are all 0,
and that explains why we should have
$(R_\sigma)_{j,k}=(R_\sigma)_{j+1,k}$ for $1\leqslant k\leqslant n$;
the same explanation can also be done for the second part of (c).\\
When the constrain (iii) appears, let us call the couple $(i,j)$ a
{\bf position} of constrain (iii).

\begin{rem}\label{rem5}
If two horizontal (resp. vertical) consecutive boxes of a matrix in
$\bR_n^2$ differ by one, then it is also the same for any
consecutive horizontal (resp. vertical) boxes above (resp. on the
right).
\end{rem}

As an immediate corollary of Proposition \ref{ens} we get

\begin{lem}\label{lemme1}
Let $\sigma,\sigma_1$ and $\sigma_2$ be involutions such that
$\sigma=\sigma_1.\sigma_2$ and $L(\sigma)=L(\sigma_1)+L(\sigma_2)$,
then $R_\sigma=R_{\sigma_1}+R_{\sigma_2}$; in particular we have
$\sigma_1,\sigma_2 \preceq \sigma$.
\end{lem}

\begin{proof}
The hypothesis $L(\sigma)=L(\sigma_1)+L(\sigma_2)$ means exactly
that any integer appearing a cycle of $\sigma_1$ does not appear in
any cycle of $\sigma_2$ and conversely [note that it is also
equivalent to say $\sigma_1.\sigma_2=\sigma_2.\sigma_1=\sigma$];
this means in particular that when the coefficient 1 appears in the
matrix $R_{\sigma_1}$ for the coordinate $(i,\sigma_1(i))$, then it
can not appear in the $i^{th}$ line and in the $\sigma_1(i)^{th}$
column of $R_{\sigma_2}$ and conversely; therefore we get
$N_\sigma=N_{\sigma_1}+N_{\sigma_2}$ and the result follows.

\end{proof}


\subsection{}\label{4.3a}

Define the following partial order on $M_{n\times n}(\Ze^+).$ For
$A,B\in M_{n\times n}(\Ze^+)$ put $A\preceq B$\index{$A\preceq B$}
if for any $i,j\ :\ 1\leqslant i,j\leqslant n$ one has
$A_{i,j}\leqslant B_{i,j}.$ \\The restriction of this order  to
$\bR_n^2$ induces a partial order on $\bS_n^2$ by setting
$\sigma'\preceq \sigma$ if $R_{\sigma'}\preceq R_{\sigma}$ for
$\sigma,\sigma'\in \bS_n^2.$ By \cite[3.5]{Mel1} this partial order describes the closures of $B.N_\sigma$ for $\sigma\in\bS_n^2$. Combining \cite[3.5]{Mel1} with Remark \ref{rem3} we get

\begin{thm} \label{thm1}
For any $\sigma\in \bS_n^2$, one has
$$\overline{B.N}_\sigma=\coprod\limits_{\sigma\pr\preceq\sigma}B.N_{\sigma\pr}.$$
In particular, for $T\in\bT_n^2$
$$\V_T=\coprod\limits_{\substack{\sigma\pr\preceq\sigma_T \\L(\sigma')=L(\sigma_T)}}B.N_{\sigma\pr}.$$
\end{thm}

\subsection{}\label{4.3b}

Let $\pi_{i,j}:\nil_{n}\rightarrow\nil_{j-i+1}$. If we denote by
$\hat{\pi}_{s,t}:\nil_{j-i+1}\rightarrow\nil_{t-s+1}$ the same
projection, but with the starting-space $\nil_{j-i+1}$, then we can
easily check the following relation:

\begin{eqnarray}\label{relation}
\hat{\pi}_{s,t}\circ\pi_{i,j}=\pi_{s+i-1,t+i-1}.
\end{eqnarray}

\noindent Now if $R\in \bR_n^2$, it is obvious by Remark \ref{ones} that
$\pi_{i,j}(R)$
fulfills the constraints (i), (ii) and (iii) of Proposition \ref{ens}.
Thus,  we get
\begin{lem}
If $R\in \bR_n^2$, then $\pi_{i,j}(R)\in \bR_{j-i+1}^2$ for
 $1\leqslant i\leqslant j\leqslant n$.
\end{lem}

\noindent Obviously, the converse is not true, as one can check for
the matrix $\tiny{\left(
             \begin{array}{ccc}
               0 & 1 & 2 \\
               0 & 0 & 1 \\
               0 & 0 & 0 \\
             \end{array}
      \right)}.$\\

      By this lemma, for any $R_\sigma\in \bR_n^2$, we have
$\pi_{i,j}(R_\sigma)\in \bR_{j-i+1}^2$; therefore $\pi_{i,j}$
induces a natural map from $\bS_n^2$ onto $\bS_{\langle
i,j\rangle}^2\cong\bS_{j-i+1}^2$, symmetric group of the set
$\{i,\cdots,j\}$. This projection will be also denoted by $\pi_{i,j}$.
Moreover by  (\ref{relation}), and Remark \ref{ones} one gets
immediately:

\begin{eqnarray}\label{relation1}
\pi_{i,j}(N_{\sigma})=N_{\pi_{i,j}(\sigma)}\ {\rm and}\
\pi_{i,j}(R_{\sigma})=R_{\pi_{i,j}(\sigma)}.
\end{eqnarray}

\noindent Note that the resulting element
$\pi_{i,j}(\sigma)$ is obtained from $\sigma$ by deleting all the cycles in which
at least one entry does not belong to $\{i,\ldots, j\}.$
 For every $\delta\in
\bS_{\langle i,j\rangle}^2$, any element $\sigma\in \pi_{i,j}^{-1}(\delta)$ will
be called a {\bf lifting} of $\delta$. In the same way we will call
the matrix $R_\sigma$ a lifting of $R_\delta$.

\begin{rem}\label{rem4}

\begin{enumerate}
\item[(i)] We will consider sometimes $\sigma\in \bS_{\langle i,j\rangle}^2$ as an element
of $\bS_n^2$ (cf. proofs of Proposition \ref{prop2}, Lemma
\ref{lemme2} and Theorem \ref{intersection});
 in particular, with the description above we have
 $\sigma=\pi_{i,j}(\sigma)$

\item[(ii)] By note (i) and Lemma \ref{lemme1}  for any
$\delta\in\bS_{\langle i,j\rangle}^2$ and any $\sigma$ its lifting in $\bS_n^2$
one has $\delta\preceq \sigma.$

\item[(iii)] By the relations
(\ref{relation1}), the projection $\pi_{i,j}$ respect the order
$\preceq$: If $\sigma_1\preceq \sigma_2$, then
$\pi_{i,j}(\sigma_1)\preceq \pi_{i,j}(\sigma_2)$.

\end{enumerate}
\end{rem}

\subsection{}\label{4.3c}

Put \index{$\bS_n^2(k)$}$\bS_n^2(k):=\{\sigma\in \bS_n^2\ |\
L(\sigma)=k\}$ and respectively $\bT_n^2(k):=\{T\in\bT_n^2\ |\
\sh(T)=(n-k,k)^*\}$. As a corollary of partial order $\preceq$ on
$\bS_n^2$ we get

\begin{prop}\label{prop2}
$\sigma_o(k):=(1,n-k+1)(2,n-k+2)\ldots(k,n)$ is the unique minimal
involution in $\bS_n^2(k)$ and for any $\sigma\in \bS_n^2(k)$ one
has $\sigma_o(k)\preceq \sigma.$ In particular, for any
$S,T\in\bT_n^2(k)$ one has $\V_T\cap\V_S\ne\emptyset.$
\end{prop}

\begin{proof}
Note that $N_{\sigma_o(k)}$ and respectively $R_{\sigma_o(k)}$ are
\begin{eqnarray}  \centering
\psset{unit=0.5cm}
\begin{pspicture}(-2,-7)(10,1)

\rput{U}(0,1){\tiny$\underbrace{n-k+1}$} \rput{U}(4,-3){\{\tiny$k$}

\rput{U}(0,0){$1$} \rput{U}(1,-1){$1$} \rput{U}(3,-3){$1$}
\psline[linestyle=dotted,dotsep=3pt]{-}(1.5,-1.5)(2.5,-2.5)

\rput{U}(-1,0){$0$} \rput{U}(0,-1){$0$} \rput{U}(3,-4){$0$}
\psline[linestyle=dotted,dotsep=3pt]{-}(0.5,-1.5)(2.5,-3.5)

\rput{U}(1,0){$0$} \rput{U}(3,-2){$0$}
\psline[linestyle=dotted,dotsep=3pt]{-}(1.5,-0.5)(2.5,-1.5)

\rput{U}(3,0){$0$}
\psline[linestyle=dotted,dotsep=3pt]{-}(1.5,0)(2.5,0)
\psline[linestyle=dotted,dotsep=3pt]{-}(3,-0.5)(3,-1.5)

\psline[linestyle=dotted,dotsep=3pt]{-}(3,-4.5)(3,-5.5)
\rput{U}(3,-6){$0$}

\rput{U}(-3,0){$0$}
\psline[linestyle=dotted,dotsep=3pt]{-}(-2.5,0)(-1.5,0)

\psline[linestyle=dotted,dotsep=3pt]{-}(-3,-0.5)(-3,-5.5)
\rput{U}(-3,-6){$0$}
\psline[linestyle=dotted,dotsep=3pt]{-}(-2.5,-6)(2.5,-6)

\rput{U}(-5,-3){$N_{\sigma_o(k)}=$}


\rput{U}(12,0){$1$} \rput{U}(13,-1){$1$} \rput{U}(15,-3){$1$}
\psline[linestyle=dotted,dotsep=3pt]{-}(13.5,-1.5)(14.5,-2.5)

\rput{U}(11,0){$0$} \rput{U}(12,-1){$0$} \rput{U}(15,-4){$0$}
\psline[linestyle=dotted,dotsep=3pt]{-}(12.5,-1.5)(14.5,-3.5)

\rput{U}(13,0){$2$} \rput{U}(15,-2){$2$}
\psline[linestyle=dotted,dotsep=3pt]{-}(13.5,-0.5)(14.5,-1.5)

\rput{U}(15,0){$k$}
\psline[linestyle=dotted,dotsep=3pt]{-}(13.5,0)(14.5,0)
\psline[linestyle=dotted,dotsep=3pt]{-}(15,-0.5)(15,-1.5)

\psline[linestyle=dotted,dotsep=3pt]{-}(15,-4.5)(15,-5.5)
\rput{U}(15,-6){$0$}

\rput{U}(9,0){$0$}
\psline[linestyle=dotted,dotsep=3pt]{-}(9.5,0)(10.5,0)

\psline[linestyle=dotted,dotsep=3pt]{-}(9,-0.5)(9,-5.5)
\rput{U}(9,-6){$0$}
\psline[linestyle=dotted,dotsep=3pt]{-}(9.5,-6)(14.5,-6)

\rput{U}(7,-3){$R_{\sigma_o(k)}=$} \rput{U}(5,-3){$,$}

\end{pspicture}
\end{eqnarray}
so that
$$(R_{\sigma_o(k)})_{i,j}=\left\{
\begin{array}{ll}
 j-i+1-(n-k) &{\rm if}\   j-i>n-k-1\\
  0 &{\rm otherwise}\end{array}\right.$$
On the other hand by Proposition \ref{ens} (ii) for any $\sigma\in
\bS_n^2$ one has $(R_\sigma)_{i,j}\geqslant
(R_\sigma)_{i-1,j}-1\geqslant (R_\sigma)_{i-2,j}-2\cdots\geqslant
(R_\sigma)_{1,j}-(i-1).$ In turn $(R_\sigma)_{1, j}\geqslant
(R_\sigma)_{1,j+1}-1\geqslant\ldots\geqslant (R_\sigma)_{1,n}-(n-j)$
so that $(R_\sigma)_{i,j}\geqslant (R_\sigma)_{1,n}-(n-j+i-1).$
Thus, for any $\sigma\in\bS_n^2(k)$ one has
$(R_\sigma)_{i,j}\geqslant j-i+1-(n-k).$ As well one has
$(R_\sigma)_{i,j}\geqslant 0$ so that $(R_\sigma)_{i,j}\geqslant
\max\{0, j-i+1-(n-k)\}=(R_{\sigma_o(k)})_{i, j}.$ Thus,
$\sigma\succeq \sigma_o(k).$

The second part is now a corollary of this result and Theorem
\ref{thm1}.

\end{proof}

\subsection{}\label{4.4}

Given $\sigma,\sigma\pr\in \bS_n^2$ we define $R_{\sigma,\sigma\pr}$
by
\begin{eqnarray}(R_{\sigma,\sigma\pr})_{i,j}:=\min\{(R_\sigma)_{i,j},(R_{\sigma\pr})_{i,j}\}.
\end{eqnarray}
One has

\begin{thm}[{\bf Main Theorem}]\label{intersection}
For any $\sigma,\sigma\pr\in \bS_n^2$ one has
$$\overline
{B.N}_\sigma\cap\overline{B.N}_{\sigma\pr}=\coprod\limits_{R_\varsigma\preceq
R_{\sigma,\sigma'}}B.N_\varsigma.$$ This intersection is
 irreducible if and only if $R_{\sigma,\sigma\pr}\in \bR_n^2.$
\end{thm}

\begin{proof}
To establish this equivalence we need only to prove the ``only if"
part and to do this we need some preliminary result.

\begin{lem}\label{lemme2}
Suppose that $\overline
{B.N}_\sigma\cap\overline{B.N}_{\sigma\pr}$ is irreducible. Denote
$B'$ the Borel subgroup in ${\rm GL}_{j-i+1}$. Then $\overline
{B'.N}_{\pi_{i,j}(\sigma)}\cap\overline{B'.N}_{\pi_{i,j}(\sigma\pr)}$
is irreducible.
\end{lem}

\begin{proof}
 Let $\alpha,\beta$ be
two maximal involutions in $\bS_{\langle i,j\rangle}^2$ such $\alpha,\beta\preceq
\pi_{i,j}(\sigma),\pi_{i,j}(\sigma')$. By Remark \ref{rem4} (ii), we
have also $\alpha,\beta \preceq \sigma, \sigma'$. By hypothesis we
have $\overline
{B.N}_\sigma\cap\overline{B.N}_{\sigma\pr}=\overline{B.N}_{\delta}$
for an element $\delta\in \bS_n^2$. In particular we get
$\alpha,\beta \preceq \delta$. By Remarks \ref{rem4} (i) and (iii)
we get $\alpha=\pi_{i,j}(\alpha), \beta=\pi_{i,j}(\beta)\preceq
\pi_{i,j}(\delta)\preceq \pi_{i,j}(\sigma),\pi_{i,j}(\sigma')$.
Since $\alpha$ and $\beta$ are maximal, we get
$\alpha=\beta=\pi_{i,j}(\delta)$.
\end{proof}
We prove the theorem by induction on $n.$
For $n=3$ all the intersections are irreducible so that
the claim is trivially true.

Let now $n$ be minimal such that $\overline
{B.N}_\sigma\cap\overline{B.N}_{\sigma\pr}$ is irreducible and
$R_{\sigma,\sigma\pr}\notin\bR_n^2$. Note that constrains (i) and
(ii) of Proposition \ref{ens} are satisfied by any
$R_{\sigma,\sigma'}.$  If $R_{\sigma,\sigma\pr}\notin\bR_n^2$ then
at least one of the conditions (a), (b) and (c) of the constrain
(iii) of Proposition \ref{ens} is not fulfilled. By symmetry around
the anti diagonal it is enough to check only Condition (a) and the
first part of Condition (c).

As for the first relation in (\ref{relation1}), we can easily check
that

\begin{eqnarray}\label{eq6}
R_{\pi_{i,j}(\sigma),\pi_{i,j}(\sigma')}=\pi_{i,j}(R_{\sigma,\sigma'}).
\end{eqnarray}

\noindent Let $B'$ be the Borel subgroup of ${\rm GL}_{n-1}.$ By
Lemma \ref{lemme2} and Relation (\ref{eq6}), we get that the
varieties $\overline
{B'.N}_{\pi_{1,n-1}(\sigma)}\cap\overline{B'.N}_{\pi_{1,n-1}(\sigma\pr)}$,
$\overline
{B'.N}_{\pi_{2,n}(\sigma)}\cap\overline{B'.N}_{\pi_{2,n}(\sigma\pr)}$
are irreducible. Thus by induction hypothesis

\begin{eqnarray}\label{eq7}
\pi_{1,n-1}(R_{\sigma,\sigma'}),\ \pi_{2,n}(R_{\sigma,\sigma'})\in
\bR_{n-1}^2,
\end{eqnarray}
Put $\zeta\in\bS_{n-1}^2$ to be such that
$R_\zeta=\pi_{1,n-1}(R_{\sigma,\sigma'})$ and $\eta\in \bS_{\langle
2,n\rangle}^2$ be such that $R_\eta=\pi_{2,n}(R_{\sigma,\sigma'}).$

Suppose that $R_{\sigma,\sigma\pr}\notin\bR_n^2$, denote $(i_o,j_o)$
the position of a constrain (iii) $\tiny{
\begin{tabular}{|c|c|} \hline
k & k+1 \\
\hline
k & k \\
\hline
\end{tabular}}$ which is not satisfied
by the matrix $R_{\sigma,\sigma\pr}$.

\vspace{0.5cm} \underline{Condition} (a): If the first part of
Condition (a) is not satisfied, it means that we can find two
horizontal consecutive boxes below of the two boxes {\tiny
\begin{tabular}{|c|c|}
        \hline
        k & k \\
        \hline
\end{tabular}}
which differ by one; but these two boxes  and {\tiny
\begin{tabular}{|c|c|}
        \hline
        k & k \\
        \hline
\end{tabular}} will lies in $\pi_{2,n}(R_{\sigma,\sigma\pr})\in\bR_{n-1}^2$, which
is impossible by Remark \ref{rem5}.

Now if the second part of Condition (a) is not satisfied, it means
that we can find two equal vertical consecutive boxes {\tiny
\begin{tabular}{|c|}
  \hline
  m \\
  \hline
  m \\
  \hline
\end{tabular}} on the right of the boxes {\tiny
\begin{tabular}{|c|}
  \hline
  k+1 \\
  \hline
  k \\
  \hline
\end{tabular}} . By Relation (\ref{eq7}), these four last boxes can not lie
inside $\pi_{1,n-1}(R_{\sigma,\sigma'}),\
\pi_{2,n}(R_{\sigma,\sigma'})$; we deduce in particular that $i_o=1$
and that the boxes {\tiny \begin{tabular}{|c|}
  \hline
  m \\
  \hline
  m \\
  \hline
\end{tabular}} belong to the last column. Since $R_{\sigma,\sigma'}$
satisfies Condition (ii) of Proposition \ref{ens}, the ``North-East"
corner of $R_{\sigma,\sigma'}$ must be {\tiny
\begin{tabular}{|c|c|}
        \hline
       m & m \\
       \hline
       m-1 & m \\
        \hline
\end{tabular}} . Now if we look at $\zeta$ (resp. $\
\eta$) as its own lifting in $\bS_n^2$, then its configuration in
the ``North-East" corner will be of the following {\tiny
\begin{tabular}{|c|c|}
        \hline
       m & m \\
       \hline
       m-1 & m-1 \\
        \hline
\end{tabular}} (resp. {\tiny
\begin{tabular}{|c|c|}
        \hline
       m-1 & m \\
       \hline
       m-1 & m \\
        \hline
\end{tabular}} ). Since the intersection is irreducible, we should
find $\delta\in\bS_n^2$ such that $\delta\succeq \zeta,\eta$ and
 $R_\delta\preceq R_{\sigma,\sigma'}.$ Since
 $(R_{\zeta})_{2,n-1}=(R_{\sigma,\sigma'})_{2,n-1}=m-1$ we get that
 also
 $(R_\delta)_{2,n-1}=m-1.$ Since
 $(R_{\zeta})_{1,n-1}=(R_{\sigma,\sigma'})_{1,n-1}=m$ we get that
 also $(R_\delta)_{1,n-1}=m.$ Since $(R_{\eta})_{2,n}=(R_{\sigma,\sigma'})_{2,n}=m$ we get that
 also $(R_\delta)_{2,n}=m.$ But then by Remark \ref{rem5} the
  ``North-East" corner of $R_\delta$ should be of the
following configuration {\tiny
\begin{tabular}{|c|c|}
        \hline
       m & m+1 \\
       \hline
       m-1 & m \\
        \hline
\end{tabular}},
this is impossible since $(R_\delta)_{1,n}\leqslant
(R_{\sigma,\sigma'})_{1,n}=m$.

\vspace{0.5cm} \underline{Condition} (c): Suppose that the first
part of Condition (c) is not satisfied, it means that we can find
two vertical consecutive boxes {\tiny \begin{tabular}{|c|}
  \hline
  m+1 \\
  \hline
  m \\
  \hline
\end{tabular}} lying in the $j_o^{th}$ and $(j_o+1)^{th}$ lines. As
above this problem can not appear inside the matrices
$\pi_{1,n-1}(R_{\sigma,\sigma'})$ and $\
\pi_{2,n}(R_{\sigma,\sigma'})$; we deduce in particular that $i_o=1$
and that the boxes {\tiny
\begin{tabular}{|c|}
  \hline
  m+1 \\
  \hline
  m \\
  \hline
\end{tabular}} lie on the last column. Since $R_{\sigma,\sigma'}$
satisfies Condition (ii) of Proposition \ref{ens}, on the right side
of the $j_o^{th}$ and $(j_o+1)^{th}$ lines of $R_{\sigma,\sigma'}$
we should find {\tiny
\begin{tabular}{|c|c|}
        \hline
       m & m+1 \\
       \hline
       m & m \\
        \hline
\end{tabular}}. Let us draw its configuration

\begin{eqnarray}\label{eq8}
\centering \psset{unit=0.5cm}
\begin{pspicture}(-5,-4)(8,0)

\rput{U}(0,0){\tiny{
\begin{tabular}{|c|c|} \hline
k & k+1 \\
\hline
k & k \\
\hline
\end{tabular}}}

\psline[linestyle=dotted,dotsep=5pt]{-}(-0.6,-0.7)(-0.6,-4)
\psline[linestyle=dotted,dotsep=5pt]{-}(0.6,-0.7)(0.6,-4)

\rput{U}(9,-3.1){{\tiny
\begin{tabular}{|c|c|}
        \hline
       m & m+1 \\
       \hline
       m & m \\
        \hline
\end{tabular}}}

\psline[linestyle=dotted,dotsep=5pt]{-}(-1,-2.9)(7.5,-2.9)
\psline[linestyle=dotted,dotsep=5pt]{-}(-1,-3.5)(7.5,-3.5)

\rput{U}(-5,-2.6){$R_{\sigma,\sigma'}$=}

\end{pspicture}
\end{eqnarray}

In the same way if we look at $R_\zeta$ and $R_\eta$ as elements of
$\bR_n^2$, then their configurations will be of the following

\begin{eqnarray}\label{eq9}
\centering \psset{unit=0.5cm}
\begin{pspicture}(-5,-4)(8,0)

\rput{U}(0,0){\tiny{
\begin{tabular}{|c|c|} \hline
k & k+1 \\
\hline
k & k \\
\hline
\end{tabular}}}

\psline[linestyle=dotted,dotsep=5pt]{-}(-0.6,-0.7)(-0.6,-4)
\psline[linestyle=dotted,dotsep=5pt]{-}(0.6,-0.7)(0.6,-4)

\rput{U}(9,-3.1){{\tiny
\begin{tabular}{|c|c|}
        \hline
       m & m \\
       \hline
       m & m \\
        \hline
\end{tabular}}}

\psline[linestyle=dotted,dotsep=5pt]{-}(-1,-2.9)(7.5,-2.9)
\psline[linestyle=dotted,dotsep=5pt]{-}(-1,-3.5)(7.5,-3.5)

\rput{U}(-5,-2.6){$R_\zeta$=}

\end{pspicture}
\end{eqnarray}

\noindent and

\begin{eqnarray}\label{eq10}
\centering \psset{unit=0.5cm}
\begin{pspicture}(-5,-4)(8,0)

\rput{U}(0,0){\tiny{
\begin{tabular}{|c|c|} \hline
k & k \\
\hline
k & k \\
\hline
\end{tabular}}}

\psline[linestyle=dotted,dotsep=5pt]{-}(-0.6,-0.7)(-0.6,-4)
\psline[linestyle=dotted,dotsep=5pt]{-}(0.6,-0.7)(0.6,-4)

\rput{U}(9,-3.1){{\tiny
\begin{tabular}{|c|c|}
        \hline
       m & m+1 \\
       \hline
       m & m \\
        \hline
\end{tabular}}}

\psline[linestyle=dotted,dotsep=5pt]{-}(-1,-2.9)(7.5,-2.9)
\psline[linestyle=dotted,dotsep=5pt]{-}(-1,-3.5)(7.5,-3.5)

\rput{U}(-5,-2.6){$R_\eta$=}

\end{pspicture}
\end{eqnarray}

\noindent Since $\delta\succeq \zeta,\eta$ and $R_\delta\preceq
R_{\sigma,\sigma'}$  combining the pictures in (\ref{eq8}),
(\ref{eq9}) and (\ref{eq10}), we get

\begin{eqnarray}
\centering \psset{unit=0.5cm}
\begin{pspicture}(-5,-4)(8,0)

\rput{U}(0,0){\tiny{
\begin{tabular}{|c|c|} \hline
k & k+1 \\
\hline
k & k \\
\hline
\end{tabular}}}

\psline[linestyle=dotted,dotsep=5pt]{-}(-0.6,-0.7)(-0.6,-4)
\psline[linestyle=dotted,dotsep=5pt]{-}(0.6,-0.7)(0.6,-4)

\rput{U}(9,-3.1){{\tiny
\begin{tabular}{|c|c|}
        \hline
       m & m+1 \\
       \hline
       m & m \\
        \hline
\end{tabular}}}

\psline[linestyle=dotted,dotsep=5pt]{-}(-1,-2.9)(7.5,-2.9)
\psline[linestyle=dotted,dotsep=5pt]{-}(-1,-3.5)(7.5,-3.5)

\rput{U}(-5,-2.6){$R_\delta$=}

\end{pspicture}
\end{eqnarray}

\noindent which is impossible, because it does not satisfy Condition
(iii) (c).

\end{proof}

\subsection{}\label{4.5}

Let us apply the previous subsection to the elements of the form
$\sigma_T$ to show that in general the intersection
$\V_T\cap\V_{T\pr}$  is reducible and
not equidimensional.

\begin{ex}\label{ex}
\begin{enumerate}
\item[(i)] For $n\leqslant 4$ all the intersections of $B$-orbit closures
of nilpotent order 2 are irreducible. The first examples of reducible
intersections of $B$ orbit closures occur  in $n=5.$
In particular there is the unique example of the reducible intersection
of orbital varieties and it is

$$T=\begin{tabular}{|c|c|}
      \hline
      1 & 2 \\
      \hline
      3 & 4 \\
      \hline
      5  \\
      \cline{1-1}
    \end{tabular}
,\ R_{\sigma_T}=\left(\tiny{\begin{array}{ccccc}
0&1&1&2&2\\
0&0&0&1&1\\
0&0&0&1&1\\
0&0&0&0&0\\
0&0&0&0&0\\\end{array}}\right)\ {\rm and } \
 T'=\begin{tabular}{|c|c|}
      \hline
      1 & 3 \\
      \hline
      2 & 5 \\
      \hline
      4  \\
      \cline{1-1}
    \end{tabular}
,\ R_{\sigma_{T'}}=\left(\tiny{\begin{array}{ccccc}
0&0&1&1&2\\
0&0&1&1&2\\
0&0&0&0&1\\
0&0&0&0&1\\
0&0&0&0&0\\\end{array}}\right) $$ so that
$$R_{\sigma_{T},\sigma_{T'}}=\left(\tiny{\begin{array}{ccccc}
0&0&1&1&2\\
0&0&0&1&1\\
0&0&0&0&1\\
0&0&0&0&0\\
0&0&0&0&0\\\end{array}}\right)$$
 Since
$(R_{\sigma_{T},\sigma_{T'}})_{1,3}=(R_{\sigma_{T},\sigma_{T'}})_{1,2}+1=(R_{\sigma_{T},\sigma_{T'}})_{2,2}+1=(R_{\sigma_{T},\sigma_{T'}})_{2,3}+1$
and
$(R_{\sigma_{T},\sigma_{T'}})_{3,5}=(R_{\sigma_{T},\sigma_{T'}})_{4,5}+1$
we get that $R_{\sigma_{T},\sigma_{T'}}$ does not satisfy condition
(iii)-(c) of Proposition \ref{ens}, therefore
$R_{\sigma_{T},\sigma_{T'}}\not\in\bR_5^2$. As well
$(R_{\sigma_T,\sigma_{T'}})_{1,4},\
(R_{\sigma_T,\sigma_{T'}})_{2,5}$ do not satisfy Remark \ref{rem5}.
Accordingly we find three maximal elements
 $R,R',R\prpr\in \bR_5^2$ for which $R,R',R\prpr\prec R_{\sigma_{T},\sigma_{T'}}$

 $$R=R_{(1,3)(2,5)}=\left(\tiny{\begin{array}{ccccc}
0&0&1&1&2\\
0&0&0&0&1\\
0&0&0&0&0\\
0&0&0&0&0\\
0&0&0&0&0\\\end{array}}\right),\quad
 R'=R_{(1,4)(3,5)}=\left(\tiny{\begin{array}{ccccc}
0&0&0&1&2\\
0&0&0&0&1\\
0&0&0&0&1\\
0&0&0&0&0\\
0&0&0&0&0\\\end{array}}\right),$$
$$
 R\prpr=R_{(1,5)(2,4)}=\left(\tiny{\begin{array}{ccccc}
0&0&0&1&2\\
0&0&0&1&1\\
0&0&0&0&0\\
0&0&0&0&0\\
0&0&0&0&0\\\end{array}}\right)$$
 Note that
 $\dim(B.N_{(1,3)(2,5)})=\dim(B.N_{(1,4)(3,5)})=\dim(B.N_{(1,5)(2,4)})=4$
 so that $\V_T\cap\V_{T'}$ contains three components of codimension 2.

\item[(ii)]
 The first example of non-equidimensional intersection of orbital varieties occurs in
 $n=6$ and it is
$$T=\begin{tabular}{|c|c|}
      \hline
      1 & 3 \\
      \hline
      2 & 6 \\
      \hline
      4  \\
      \cline{1-1}
      5  \\
      \cline{1-1}

    \end{tabular}
\ ,\quad R_{\sigma_{T}}=\left(\tiny{\begin{array}{cccccc}
0&0&1&1&1&2\\
0&0&1&1&1&2\\
0&0&0&0&0&1\\
0&0&0&0&0&1\\
0&0&0&0&0&1\\
0&0&0&0&0&0\\
\end{array}}\right)$$  and

$$
 T'=\begin{tabular}{|c|c|}
      \hline
      1 & 2 \\
      \hline
      3 & 5 \\
      \hline
      4  \\
      \cline{1-1}
      6  \\
      \cline{1-1}

    \end{tabular}
\ ,\quad R_{\sigma_{T'}}=\left(\tiny{\begin{array}{cccccc}
0&0&1&1&2&2\\
0&0&0&0&1&1\\
0&0&0&0&1&1\\
0&0&0&0&1&1\\
0&0&0&0&0&0\\
0&0&0&0&0&0\\
\end{array}}\right) $$ so that
$$R_{\sigma_{T},\sigma_{T'}}=\left(\tiny{\begin{array}{cccccc}
0&0&1&1&1&2\\
0&0&0&0&1&1\\
0&0&0&0&0&1\\
0&0&0&0&0&1\\
0&0&0&0&0&0\\
0&0&0&0&0&0\\\end{array}}\right)$$

\noindent Since
$(R_{\sigma_{T},\sigma_{T'}})_{1,3}=(R_{\sigma_{T},\sigma_{T'}})_{1,2}+1=(R_{\sigma_{T},\sigma_{T'}})_{2,2}+1=(R_{\sigma_{T},\sigma_{T'}})_{2,3}+1$
and
$(R_{\sigma_{T},\sigma_{T'}})_{1,5}=(R_{\sigma_{T},\sigma_{T'}})_{2,5}$
we get that $R_{\sigma_{T},\sigma_{T'}}$ does not satisfy condition
(iii) (a) of Proposition \ref{ens}  and
$(R_{\sigma_T,\sigma_{T'}})_{1,5},\
(R_{\sigma_T,\sigma_{T'}})_{2,6}$ do not satisfy  Remark \ref{rem5}
so that $R_{\sigma_{T},\sigma_{T'}}\not\in\bR_6^2$ and the maximal
elements
 $R,R'\in \bR_6^2$ for which $R,R'\prec R_{\sigma_{T},\sigma_{T'}}$ are
 $$R=R_{(1,3)(4,6)}=\left(\tiny{\begin{array}{cccccc}
0&0&1&1&1&2\\
0&0&0&0&0&1\\
0&0&0&0&0&1\\
0&0&0&0&0&1\\
0&0&0&0&0&0\\
0&0&0&0&0&0\\\end{array}}\right)$$ and

$$R'=R_{(1,6)(2,5)}=\left(\tiny{\begin{array}{cccccc}
0&0&0&0&1&2\\
0&0&0&0&1&1\\
0&0&0&0&0&0\\
0&0&0&0&0&0\\
0&0&0&0&0&0\\
0&0&0&0&0&0\\
\end{array}}\right).$$
Note that
 $\dim(B.N_{(1,3)(4,6)})=6$ and $\dim(B.N_{(1,6)(2,5)})=4$
 so that $\V_T\cap\V_{T'}$ contains one component of codimension 2
 and another component of codimension 4.
\end{enumerate}

\end{ex}

\section{Some other counter-examples}
\subsection{Cell graphs }\label{4.6}

Let $T\in \bT_{\lambda} $ be a standard tableau and $\mathcal{C}_T$ its corresponding left cell (cf. Section \ref{1.4}).
 Steinberg's construction provides the way to construct  $\mathcal{V}_T$ with the help of elements of $\mathcal{C}_T$.
In \cite{M-P}, we got another geometric interpretation of $C_T$:

\begin{thm}\label{thm2}(\cite{M-P}) -
Let $T\in \bT_{\lambda}$ and let $w={\rm RS}(T,T\pr)\in \mathcal{C}_T$. Then for
a  $x\in\mathcal{V}_{T}\cap B.(\nil\cap^w\nil)$ in general position, the unique Schubert cell whose intersection with
the irreducible component $\F_{T\pr}$ of the Springer fiber is open and dense in $\F_{T\pr}$ is indexed by $w$.
\end{thm}

 The cell $\mathcal{C}_T$ in $\bS_n$ can be visualized as a  {\bf cell graph} $\Gamma_T$ where the vertices are labeled by $\bT_{\lambda}$, and
two vertices $T'$ and $T''$ are joined by an edge labeled by $k$ if  $s_k{\rm RS}(T,T\pr)={\rm RS}(T,T'')$. One can easily see ( cf. \cite{M-P}, for example) that if $T'$ and $T''$ are joined in $\Gamma_T$, then $\codim_{\F_{T'}}\F_{T\pr}\cap\F_{T''}=1.$

 Note that $T'$ and $T''$ can be joined by an edge in $\Gamma_T$
and not joined by an edge in $\Gamma_S$ for some $S,T\in\bT_\lambda$.
Is it true that $\codim_{\F_{T'}}\F_{T'}\cap\F_{T''}=1$ if and only if there exists $T\in\bT_\lambda$ such that $T'$ and $T''$ are joined by an edge in $\Gamma_T$?

The answer is negative as we show by the example below.

 As we show in \cite{M-P1} if $k\leqslant 2$ then $\codim_{\V_T}(\V_T\cap\V_S)=1$
 if and only if there exists $P\in\bT_{(n-k,k)^*}$ such that  $T$ and $S$
 are joined by an  edge in  $\Gamma_P$  so that the
 first example occurs in $n=6$ for $\bT_{(3,3)^*}.$ In that case $(3,3)^*=(2,2,2)$ and the corresponding orbital varieties
 are $9$-dimensional. Let us put
 $$T_{1}=\begin{tabular}{|c|c|}
           \hline
           1 & 4 \\
           \hline
           2 & 5 \\
           \hline
           3 & 6 \\
           \hline
         \end{tabular}
 ,\quad
T_{2}=\begin{tabular}{|c|c|}
           \hline
           1 & 3 \\
           \hline
           2 & 5 \\
           \hline
           4 & 6 \\
           \hline
         \end{tabular}\ ,\quad T_{3}=\begin{tabular}{|c|c|}
           \hline
           1 & 2 \\
           \hline
           3 & 5 \\
           \hline
           4 & 6 \\
           \hline
         \end{tabular}\ ,\quad
T_{4}=\begin{tabular}{|c|c|}
           \hline
           1 & 3 \\
           \hline
           2 & 4 \\
           \hline
           5 & 6 \\
           \hline
         \end{tabular}\ ,\quad T_{5}=\begin{tabular}{|c|c|}
           \hline
           1 & 2 \\
           \hline
           3 & 4 \\
           \hline
           5 & 6 \\
           \hline
         \end{tabular}\ .$$
\noindent One can check that all the  cell graphs  are the same this graph  is
\begin{center}
\begin{picture}(70,120)
\put(15,50){$T_{4}$} \put(50,20){$T_{5}$} \put(50,110){$T_{1}$}
\put(50,80){$T_{2}$} \put(80,50){$T_{3}$}
 \put(52,90){\line(0,1){17}}\put(55,97){$\scriptstyle 3$}
 \put(28,55){\line(1,1){22}}\put(30,65){$\scriptstyle 4$}
 \put(76,55){\line(-1,1){22}}\put(74,65){$\scriptstyle 2$}
 \put(48,30){\line(-1,1){19}}\put(30,38){$\scriptstyle 2$}
 \put(55,30){\line(1,1){20}}\put(71,38){$\scriptstyle 4$}
 \end{picture}
 \end{center}

\noindent On the other hand one has

$$R_{\sigma_{T_1},\sigma_{T_5}}=\left(\tiny{\begin{array}{cccccc}
0&0&0&1&2&3\\
0&0&0&1&1&2\\
0&0&0&1&1&1\\
0&0&0&0&0&0\\
0&0&0&0&0&0\\
0&0&0&0&0&0\\
\end{array}}\right)=R_{(1,5)(2,6)(3,4)}$$

\noindent and $\dim(B.N_{(1,5)(2,6)(3,4)})=8$, so that
$\codim_{\V_{T_1}}(\V_{T_1}\cap\V_{T_5})=1.$ As well the straight
computations show that
$\dim(\V_{T_1}\cap\V_{T_4})=\dim(\V_{T_1}\cap\V_{T_3})=
\dim(\V_{T_2}\cap\V_{T_5})=\dim(\V_{T_3}\cap\V_{T_4})=7$ so that all
these intersections are of codimension 2. Further,
$\V_{T_1}\cap\V_{T_4},\ \V_{T_1}\cap\V_{T_3}$ and
$\V_{T_3}\cap\V_{T_4}$ are irreducible. $\V_{T_2}\cap\V_{T_5}$ has
three components with the following dense $B$-orbits:
$B.N_{(1,3)(2,5)(4,6)},\ B.N_{(1,5)(2,4)(3,6)},$ and
$B.N_{(1,4)(2,6)(3,5)}.$
 Below we draw the graph where two vertices are joined if the corresponding intersection is of codimension 1.
\begin{figure}[htbp]
  \centering
\psset{unit=0.5cm}
\begin{tabular}{cc}
\begin{picture}(100,90)

 \put(0,10){$T_{4}$}
 \put(90,10){$T_{5}$}
 \put(0,50){$T_{2}$}
 \put(90,50){$T_{3}$}
 \put(39,85){$T_{1}$}

 \put(2,20){\line(0,1){25}}
\put(10,12){\line(1,0){75}}
\put(92,20){\line(0,1){25}}
\put(10,52){\line(1,0){75}}
\put(9,58){\line(5,4){28}}
 \put(90,20){\line(-2,3){40}}

\end{picture}
\end{tabular}
\end{figure}

\subsection{Orbital variety's construction}

Let us go back to Steinberg's construction of an orbital variety (see Section \ref{1.4}). Given $T\in\bT_\lambda$ one has
 $\mathcal{V}_{T}=\overline{B.(\nil\cap {^w \nil})}\cap {\Or}_{\lambda}$
for any $w \in \mathcal{C}_T$. Obviously,
$${\rm dim}(\overline{B.(\nil\cap {^w \nil})}\cap {\Or}_{\lambda})={\rm dim}(B.(\nil\cap {^w \nil})\cap {\Or}_{\lambda}),$$
\noindent
so that ${\rm dim}(B.(\nil\cap {^w \nil})\cap {\Or}_{\lambda})={\rm dim}(\Or_{\lambda}\cap \nil)$, therefore $B.(\nil\cap {^w \nil})\cap {\Or}_{\lambda}$ is also irreducible in $\Or_{\lambda}\cap \nil$; in particular $B.(\nil\cap {^w \nil})\cap {\Or}_{\lambda}$ is an orbital variety if and only if $B.(\nil\cap {^w \nil})\cap {\Or}_{\lambda}$ is closed in $\Or_{\lambda}\cap \nil$.
The natural questions connected to the construction
are \\
 {\bf Q1}: May be one can always find $w\in{\mathcal C}_T$ such that $\mathcal{V}_T=B.(\nil\cap {^w \nil})\cap {\Or}_{\lambda}$?\\
 {\bf Q2}: Or may be  $\mathcal{V}_T=\bigcup\limits_{y\in \mathcal{C}_T}B.(\nil\cap {^y \nil})\cap {\Or}_{\lambda}$?\\

The answers to both these questions are negative as we show by the following counter-example.

\begin{ex} Let $T=\begin{tabular}{|c|c|}
        \hline
        1 & 3  \\
        \hline
        2  \\
        \cline{1-1}
        4   \\
        \cline{1-1}
      \end{tabular}$. The corresponding left cell is given by $$\mathcal{C}_T=\{w_1=\tiny{\left(
                                          \begin{array}{cccc}
                                            1 & 2 & 3 & 4 \\
                                            4 & 2 & 3 & 1 \\
                                          \end{array}
                                        \right)},\ w_2=\tiny{\left(
                                          \begin{array}{cccc}
                                            1 & 2 & 3 & 4 \\
                                            2 & 4 & 3 & 1 \\
                                          \end{array}
                                        \right)},\ w_3=\tiny{\left(
                                          \begin{array}{cccc}
                                            1 & 2 & 3 & 4 \\
                                            4 & 2 & 1 & 3 \\
                                          \end{array}
                                        \right)}\}.$$

We draw here in green the corresponding space $\nil\cap {^w \nil}$:

\begin{figure}[h]
\centerline{\includegraphics[width=10cm, height=2.3cm]{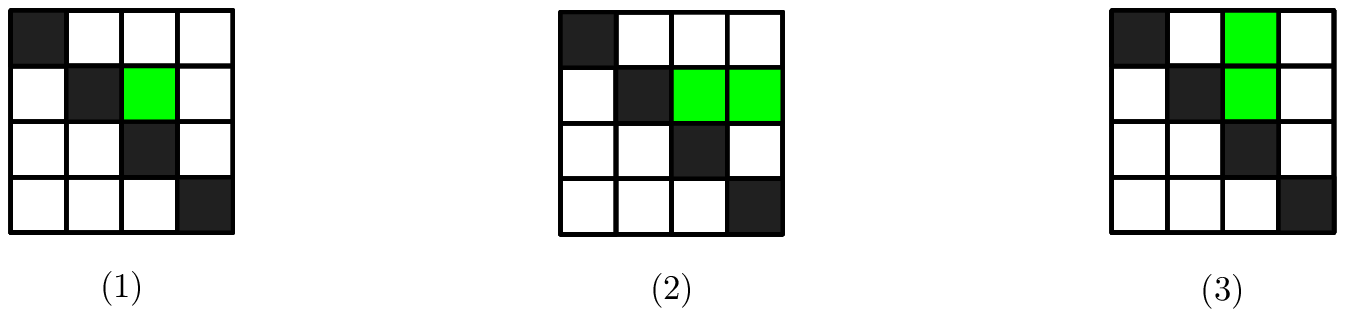}}\caption{}
\end{figure}

\noindent On the other hand by Theorem \ref{thm1},
$\mathcal{V}_T=B.N_{(2,3)}\cup B.N_{(2,4)}\cup B.N_{(1,3)}\cup B.N_{(1,4)}$.  As one can see from the picture  $N_{(1,4)}\not\in B.(\nil\cap {^w \nil})$ for $w\in \{w_1,w_2,w_3\}$.

\end{ex}

\vskip 0.2 cm {\bf Acknowledgements.}  The second author would like
to express his gratitude to A. Joseph, L\^{e} D\~{u}ng Tr\'{a}ng, H. Esnault and E. Viehweg
for the invitation to the Weizmann institute of Science, the
Abdus Salam International Centre for Theoretical Physics and the University of Duisburg-Essen where this
work was done. He would  also like to thank these institutions for their
kind hospitality and support.

\bigskip

\centerline{ INDEX OF NOTATION}
\medskip

\noindent Symbols appearing frequently are given below in order of
appearance.

\medskip

\ref{1.1} $\nil, \G_\alpha,\ \alpha_i,\ \Pi,\ \alpha_{i,j},\ B,\
    \bS_n,\ s_i,\ g.u,\ \F_x,\ \Or_x;$

\ref{1.2} $J(x),\ \Or_\lambda,\ \sh(\Or),\ \sh(T),\ \bT_\lambda,\
\F_T,\ \V_T;$

\ref{2.3} $\pi_{i,j}:\nil_n\rar\nil_{j+1-i}$

\ref{4.1} $\X_2,\ \bS_n^2,\ N_\sigma,\ \bT_n^2,\ \sigma_T;$

\ref{4.2} $L(\sigma),\ \Or_\sigma;$

\ref{4.3} $R_\sigma,\ \bR_n^2;$

\ref{4.3b} $\bS_{\langle i,j\rangle}^2,\
                 \pi_{i,j}:\bS_n^2\rightarrow\bS_{\langle i,j\rangle}^2$

\ref{4.3c} $\bS_n^2(k),\ \bT_n^2(k)$



\end{document}